# Separating Complexity Classes using Autoreducibility

Harry Buhrman[*]  Lance Fortnow[†]  Dieter van Melkebeek[‡]
CWI  The University of Chicago  The University of Chicago

Leen Torenvliet[§]
University of Amsterdam

**Abstract**

A set is autoreducible if it can be reduced to itself by a Turing machine that does not ask its own input to the oracle. We use autoreducibility to separate the polynomial-time hierarchy from polynomial space by showing that all Turing-complete sets for certain levels of the exponential-time hierarchy are autoreducible but there exists some Turing-complete set for doubly exponential space that is not.

Although we already knew how to separate these classes using diagonalization, our proofs separate classes solely by showing they have different structural properties, thus applying Post's Program to complexity theory. We feel such techniques may prove unknown separations in the future. In particular, if we could settle the question as to whether all Turing-complete sets for doubly exponential time are autoreducible, we would separate either polynomial time from polynomial space, and nondeterministic logarithmic space from nondeterministic polynomial time, or else the polynomial-time hierarchy from exponential time.

We also look at the autoreducibility of complete sets under nonadaptive, bounded query, probabilistic and nonuniform reductions. We show how settling some of these autoreducibility questions will also lead to new complexity class separations.

[*]Part of this research was done while visiting the Univ. Politècnica de Catalunya in Barcelona. Partially supported by the Dutch foundation for scientific research (NWO) through NFI Project ALADDIN, under contract NF 62-376 and a TALENT stipend. E-mail: buhrman@cwi.nl. URL: http://www.cwi.nl/cwi/people/Harry.Buhrman.html.

[†]Supported in part by NSF grant CCR 92-53582. Research partly done while visiting CWI. Email: fortnow@cs.uchicago.edu. URL: http://www.cs.uchicago.edu/~fortnow.

[‡]Partially supported by the European Union through Marie Curie Research Training Grant ERB-4001-GT-96-0783, by the U.S. National Science Foundation through grant CCR 92-53582, and by the Fields Institute. Research partly done while visiting CWI and the University of Amsterdam. Email: dieter@cs.uchicago.edu. URL: http://www.cs.uchicago.edu/~dieter.

[§]Partially supported by HC&M grant ERB4050PL93-0516. Email: leen@wins.uva.nl. URL: http://turing.wins.uva.nl/~leen.

# 1 Introduction

While complexity theorists have made great strides in understanding the structure of complexity classes, they have not yet found the proper tools to do nontrivial separation of complexity classes such as P and NP. They have developed sophisticated diagonalization, combinatorial and algebraic techniques but none of these ideas have yet proven very useful in the separation task.

Back in the early days of computability theory, Post [13] wanted to show that the set of noncomputable computably enumerable sets strictly contains the Turing-complete computably enumerable sets. In what we now call "Post's Program" (see [11, 15]), Post tried to show these classes differ by finding a property that holds for all sets in the first class but not for some set in the second.

We would like to resurrect Post's Program for separating classes in complexity theory. In particular we will show how some classes differ by showing that their complete sets have different structure. While we do not separate any classes not already separable by known diagonalization techniques, we feel that refinements to our techniques may yield some new separation results.

In this paper we will concentrate on the property known as "autoreducibility." A set $A$ is autoreducible if we can decide whether an input $x$ belongs to $A$ in polynomial-time by making queries about membership of strings different from $x$ to $A$.

Trakhtenbrot [16] first looked at autoreducibility in both the computability theory and space-bounded models. Ladner [10] showed that there exist Turing-complete computably enumerable sets that are not autoreducible. Ambos-Spies [1] first transferred the notion of autoreducibility to the polynomial-time setting. More recently, Yao [18] and Beigel and Feigenbaum [5] have studied a probabilistic variant of autoreducibility known as "coherence."

In this paper, we ask for what complexity classes do all the complete sets have the autoreducibility property. In particular we show:

- All Turing-complete sets for $\Delta_k^{\mathrm{EXP}}$ are autoreducible for any constant $k$, where $\Delta_{k+1}^{\mathrm{EXP}}$ denotes the sets that are exponential-time Turing-reducible to $\Sigma_k^{\mathrm{P}}$.

- There exists a Turing-complete set for doubly exponential space that is not autoreducible.

Since the union of all sets $\Delta_k^{\mathrm{EXP}}$ coincides with exponential-time hierarchy, we obtain a separation of the exponential-time hierarchy from doubly exponential space and thus of the polynomial-time hierarchy from exponential space. Although these results also follow from the space hierarchy theorems [9] which we have known for a long time, our proof does not directly use diagonalization, rather separates the classes by showing that they have different structural properties.

Issues of relativization do not apply to this work because of oracle access (see [8]): A polynomial-time autoreduction can not view as much of the oracle as an exponential or doubly exponential computation. To illustrate this point we show that there exists an oracle relative to which some complete set for exponential time is not autoreducible.

Note that if we can settle whether the Turing-complete sets for doubly exponential time are all autoreducible one way or the other, we will have a major separation result. If there exists a Turing-complete set for doubly exponential time that is not autoreducible, then we get that the exponential-time hierarchy is strictly contained in doubly exponential time thus that the polynomial-time hierarchy is strictly contained in exponential time. If all of the Turing-complete sets for doubly exponential time are autoreducible, we get that doubly exponential time is strictly contained in doubly exponential space, and thus polynomial time strictly in polynomial space. We will also show that this assumption implies a separation of nondeterministic logarithmic space from nondeterministic polynomial time. Similar implications hold for space bounded classes (see Section 5). Autoreducibility questions about doubly exponential time and exponential space thus remain an exciting line of research.

We also study the nonadaptive variant of the problem. Our main results scale down one exponential as follows:

- All truth-table-complete sets $\Delta_k^P$ are truth-table-autoreducible for any constant $k$, where $\Delta_{k+1}^P$ denotes that sets polynomial-time Turing-reducible to $\Sigma_k^P$.

- There exists a truth-table-complete set for exponential space that is not truth-table-autoreducible.

Again, finding out whether all truth-table-complete sets for intermediate classes, namely polynomial space and exponential time, are truth-table-autoreducible, would have major implications.

In contrast to the above results we exhibit the limitations of our approach: For the restricted reducibility where we are only allowed to ask two nonadaptive queries, all complete sets for EXP, EXPSPACE, EEXP, EEXPSPACE, etc., are autoreducible.

We also argue that uniformity is crucial for our technique of separating complexity classes, because our nonautoreducibility results fail in the nonuniform setting. Razborov and Rudich [14] show that if strong pseudo-random generators exist, "natural proofs" cannot separate certain nonuniform complexity classes. Since this paper relies on uniformity in an essential way, their result does not apply.

Regarding the probabilistic variant of autoreducibility mentioned above, we can strengthen our results and construct a Turing-complete set for doubly exponential space that is not even probabilistically autoreducible. We leave the analogue of this theorem in the nonadaptive setting open: Does there exist a truth-table complete set for exponential space that is not probabilistically truth-table autoreducible? We do show that every truth-table complete set for exponential time is probabilistically truth-table autoreducible. So, a positive answer to the open question would establish that exponential time is strictly contained in exponential space. A negative answer, on the other hand, would imply a separation of nondeterministic logarithmic space from nondeterministic polynomial time.

Here is the outline of the paper: First, we introduce our notation and state some preliminaries in Section 2. Next, in Section 3 we establish our negative autoreducibility results, for the adaptive as well as the nonadaptive case. Then we prove the positive results in Section 4, where we also briefly look at the randomized and nonuniform settings. Section 5 discusses the separations that follow from our results and would follow from improvements on them. Finally, we conclude in Section 6 and mention some possible directions for further research.

## 1.1 Errata in conference version

An previous version of this paper [6] erroneously claimed proofs showing all Turing complete sets for EXPSPACE are autoreducible and all truth-table complete sets for PSPACE are nonadaptively autoreducible. Combined with the additional results in this version, we would have a separation of NL and NP (see Section 5).

However the proofs in the earlier version failed to account for the growth of the running time when recursively computing previous players' moves. We use the proof technique in Section 3 though unfortunately we get considerably weaker theorems. The original results claimed in the previous version remain important open questions as resolving them either way will yield new separation results.

## 2 Notation and Preliminaries

Most of our complexity theoretic notation is standard. We refer the reader to the textbooks by Balcázar, Díaz and Gabarró [4, 3], and by Papadimitriou [12].

We use the binary alphabet $\Sigma = \{0, 1\}$. We denote the difference of a set $A$ with a set $B$, i.e., the subset of elements of $A$ that do not belong to $B$, by $A \setminus B$.

For any integer $k \geqslant 0$, a $\Sigma_k$-formula is a Boolean expression of the form

$$\exists y_1 \in \Sigma^{n_1}, \forall y_2 \in \Sigma^{n_2}, \ldots, Q_k y_k \in \Sigma^{n_k} : \phi(y_1, y_2, \ldots, y_k, z), \tag{1}$$

where $\phi$ is a Boolean formula, $Q_i$ denotes $\exists$ if $i$ is odd, and $\forall$ otherwise, and the $n_i$'s are positive integers. We say that (1) has $k-1$ alternations. A $\Pi_k$-formula is just like (1) except starts with a $\forall$-quantifier. It also has $k-1$ alternations. A $\mathrm{QBF}_k$-formula is a $\Sigma_k$-formula (1) without free variables $z$.

For any integer $k \geqslant 0$, $\Sigma_k^\mathrm{P}$ denotes the $k$-th $\Sigma$-level of the polynomial-time hierarchy. We define these levels recursively by $\Sigma_0^\mathrm{P} = \mathrm{P}$, and $\Sigma_{k+1}^\mathrm{P} = \mathrm{NP}^{\Sigma_k^\mathrm{P}}$. The $\Delta$-levels of the polynomial-time and exponential-time hierarchy are defined as $\Delta_{k+1}^\mathrm{P} = \mathrm{P}^{\Sigma_k^\mathrm{P}}$ respectively $\Delta_{k+1}^\mathrm{EXP} = \mathrm{EXP}^{\Sigma_k^\mathrm{P}}$. The polynomial-time hierarchy PH equals the union of all sets $\Delta_k^\mathrm{P}$, and the exponential-time hierarchy EXPH similarly the union of all sets $\Delta_k^\mathrm{EXP}$.

A *reduction* of a set $A$ to a set $B$ is a polynomial-time oracle Turing machine $M$ such that $M^B = A$. We say that $A$ reduces to $B$ and write $A \leqslant_\mathrm{T}^\mathrm{P} B$ ("T" for Turing). The reduction $M$ is *nonadaptive* if the oracle queries $M$ makes on any input are independent of the oracle, i.e., the queries do not depend upon the answers to previous queries. In that case we write $A \leqslant_\mathrm{tt}^\mathrm{P} B$ ("tt" for truth-table). Reductions of functions to sets are defined similarly. If the number of queries on an input of length $n$ is bounded by $q(n)$, we write $A \leqslant_{q(n)-\mathrm{T}}^\mathrm{P} B$ respectively $A \leqslant_{q(n)-\mathrm{tt}}^\mathrm{P} B$; if it is bounded by some constant, we write $A \leqslant_\mathrm{btt}^\mathrm{P} B$ ("b" for bounded). We denote the set of queries of $M$ on input $x$ with oracle $B$ by $Q_{M^B}(x)$; in case of nonadaptive reductions, we omit the oracle $B$ in the notation. If the reduction asks only one query and answers the answer to that query, we write $A \leqslant_\mathrm{m}^\mathrm{P} B$ ("m" for many-one).

For any reducibility $\leqslant_r^\mathrm{P}$ and any complexity class $\mathcal{C}$, a set $C$ is $\leqslant_r^\mathrm{P}$-*hard* for $\mathcal{C}$ if we can $\leqslant_r^\mathrm{P}$-reduce every set $A \in \mathcal{C}$ to $C$. If in addition $C \in \mathcal{C}$, we call $C$ $\leqslant_r^\mathrm{P}$-*complete* for $\mathcal{C}$. For any integer $k \geqslant 0$, the set $\mathrm{TQBF}_k$ of all true $\mathrm{QBF}_k$-formulae is $\leqslant_\mathrm{m}^\mathrm{P}$-complete for $\Sigma_k^\mathrm{P}$. For $k = 1$, this reduces to the fact that the set SAT of satisfiable Boolean formulae is $\leqslant_\mathrm{m}^\mathrm{P}$-complete for NP.

Now we get to the key concept of this paper:

**Definition 2.1** *A set $A$ is* autoreducible *if there is a reduction $M$ of $A$ to itself that never queries its own input, i.e., for any input $x$ and any oracle $B$, $x \notin Q_{M^B}(x)$. We call such $M$ an* autoreduction *of $A$.*

We will also discuss randomized and nonuniform variants. A set is *probabilistically autoreducible* if it has a probabilistic autoreduction with bounded two-sided error. Yao [18] first studied this concept under the name "coherence". A set is *nonuniformly autoreducible* if it has an autoreduction that uses polynomial advice. For all these notions, we can consider both the adaptive and the nonadaptive case. For randomized autoreducibility, nonadaptiveness means that the queries only depend on the input and the random seed.

## 3 Nonautoreducibility Results

In this section, we show that large complexity classes have complete sets that are not autoreducible.

**Theorem 3.1** *There is a $\leqslant^{\mathrm{P}}_{2-\mathrm{T}}$-complete set for EEXPSPACE that is not autoreducible.*

Most natural classes containing EEXPSPACE, e.g., EEEXPTIME and EEEXPSPACE, also have this property.

We can even construct the complete set in Theorem 3.1 to defeat every probabilistic autoreduction:

**Theorem 3.2** *There is a $\leqslant^{\mathrm{P}}_{2-\mathrm{T}}$-complete set for EEXPSPACE that is not probabilistically autoreducible.*

In the nonadaptive setting, we obtain:

**Theorem 3.3** *There is a $\leqslant^{\mathrm{P}}_{3-\mathrm{tt}}$-complete set for EXPSPACE that is not nonadaptively autoreducible.*

Unlike in case of Theorem 3.1, our construction does not seem to yield a truth-table complete set that is not probabilistically nonadaptively autoreducible. In fact, as we shall show in Section 4.3, such a result would separate EXP from EXPSPACE. See also Section 5.

We will detail in Section 4.3 that our nonautoreducibility results do not hold in the nonuniform setting.

### 3.1 Adaptive Autoreductions

Suppose we want to construct a $\leqslant^{\mathrm{P}}_{\mathrm{T}}$-complete set $A$ for a complexity class $\mathcal{C}$ that is not autoreducible. If $\mathcal{C}$ has a $\leqslant^{\mathrm{P}}_{\mathrm{m}}$-complete set $K$, we could try to encode $K$ in $A$, and at the same time diagonalize against all autoreductions. A straightforward implementation would be to encode $K(y)$ as $A(\langle 0, y \rangle)$, and stage-wise diagonalize against all $\leqslant^{\mathrm{P}}_{\mathrm{T}}$-reductions $M$ by picking for each $M$ an input $x$ not of the form $\langle 0, y \rangle$ that is not queried during previous stages, and setting $A(x) = 1 - M^A(x)$. However, the computation of $M^A(x)$ might require deciding $K(y)$ on inputs $y$ of order of size the running time of $M$ on input $x$. Since we have to do this for all potential autoreductions $M$, we can only bound the time complexity of $A$ by a function in $t(n^{\omega(1)})$, where $t(n)$ denotes the running time of some deterministic Turing machine accepting $K$. That does not suffice to keep $A$ inside $\mathcal{C}$.

In order to solve this problem, we consider two possible coding regions at every stage: the left region $L$, containing strings of the form $\langle 0, y \rangle$, and the right region $R$, containing strings of the form $\langle 1, y \rangle$. Now, suppose we want to diagonalize against $\leqslant^{\mathrm{P}}_{\mathrm{T}}$-reduction $M$ on input $x$ (not of the form $\langle 0, y \rangle$ or $\langle 1, y \rangle$).

**Observation 3.1** *Either it is the case that for any setting of $L$ there is a setting of $R$ such that $M^A(x)$ accepts, or for any setting of $R$ there is a setting of $L$ such that $M^A(x)$ rejects.*

In the former case, we will set $A(x) = 0$ and encode $K$ in $L$ (at that stage); otherwise we will set $A(x) = 1$ and encode $K$ in $R$. In any case,

$$K(y) = A(\langle A(x), y \rangle),$$

so deciding $K$ is still easy when given $A$. Moreover, the diagonalization, i.e., determining $A(x)$, no longer requires computing $K(y)$ for large inputs $y$: We just have to decide which of the above cases holds, and that is independent of the part of $K$ we want to encode. Provided $\mathcal{C}$ is powerful enough, we can do it in $\mathcal{C}$. The price we have to pay, is that in order to force $M^A(x)$ to $1 - A(x)$, we have to set the bits in the non-coding region so as to satisfy the underlying property of Observation 3.1. In addition to determining the coding region, this still requires the knowledge of $K(y)$ for possibly large inputs $y$. We will use a slightly stronger version of Observation 3.1 to circumvent the need to compute $K(y)$ for too large inputs $y$, by grouping the quantifiers corresponding to inputs of about the same length in Observation 3.1 and rearranging them.

This is what happens in the next lemma, which we prove in a more general form, because we will need the generalization later on in Section 5.

**Lemma 3.4** *Fix a set $K$, and suppose we can decide it simultaneously in time $t(n)$ and space $(n)$. Let $\beta : \mathbb{N} \to (0, \infty)$ be a constructible monotone unbounded function, and suppose there is a deterministic Turing machine accepting TQBF that takes time $t'(n)$ and space $'(n)$ on QBF-formulae of size $2^{n^{\beta(n)}}$ with at most $\log \beta(n)$ alternations. Then there is a set $A$ such that:*

1. *$A$ is not autoreducible.*

2. *$K \leqslant^{\mathrm{P}}_{2-\mathrm{T}} A$.*

3. *We can decide $A$ simultaneously in time $O(2^{n^2} \cdot (t(n^2) + t'(n)))$ and space $O(2^{n^2} \cdot (\ ^2(n) + \ '(n)))$.*

**Proof** (of Lemma 3.4)
Fix a function $\beta$ satisfying the hypotheses of the Lemma. Let $M_1, M_2, \ldots$ be a standard enumeration of autoreductions clocked such that $M_i$ runs in time $n^{\beta(i)}$ on inputs of length $n$. Our construction starts out with $A$ being the empty set, and then adds strings to $A$ in subsequent stages $i = 1, 2, 3, \ldots$ defined by the following sequence:

$$\begin{cases} n_0 &= 0 \\ n_{i+1} &= n_i^{\beta(n_i)} + 1. \end{cases}$$

Fix and integer $i \geqslant 1$ and let $m = n_i$. For any integer $j$ such that $0 \leqslant j \leqslant \log \beta(m)$, let $I_j$ denote the set of all strings with lengths in the interval $[m, \min(m^{2^{j+1}}, m^{\beta(m)} + 1))$. Note that $\{I_j\}_{j=0}^{\log \beta(m)}$ forms a partition of the set $I$ of strings with lengths in $[m, m^{\beta(m)} + 1) = [n_i, n_{i+1})$.

During the $i$-th stage of the construction, we will encode the restriction $K|_I$ of $K$ to $I$ into $\{\langle b, y \rangle \mid b \in \Sigma \text{ and } y \in I\}$, and use the string $0^m$ for diagonalizing against $M_i$, applying the next strengthening of Observation 3.1 to do so:

**Claim 3.1** *For any set $A$, at least one of the following holds:*

$$(\forall \ell_y)_{y \in I_0}, (\exists r_y)_{y \in I_0}, (\forall \ell_y)_{y \in I_1}, (\exists r_y)_{y \in I_1},$$
$$\ldots, (\forall \ell_y)_{y \in I_{\log \beta(m)}}, (\exists r_y)_{y \in I_{\log \beta(m)}} : M_i^{A'}(0^m) \text{ accepts} \qquad (2)$$

*or*

$$(\forall r_y)_{y \in I_0}, (\exists \ell_y)_{y \in I_0}, (\forall r_y)_{y \in I_1}, (\exists \ell_y)_{y \in I_1},$$
$$\ldots, (\forall r_y)_{y \in I_{\log \beta(m)}}, (\exists \ell_y)_{y \in I_{\log \beta(m)}} : M_i^{A'}(0^m) \text{ rejects}, \qquad (3)$$

*where $A'$ denotes $A \cup \{\langle 0, y \rangle \mid y \in I \text{ and } \ell_y = 1\} \cup \{\langle 1, y \rangle \mid y \in I \text{ and } r_y = 1\}$.*

Here we use $(Q\, z_y)_{y \in Y}$ as a shorthand for $Q\, z_{y_1}, Q\, z_{y_2}, \ldots, Q\, z_{y_{|Y|}}$, where $Y = \{y_1, y_2, \ldots, y_{|Y|}\}$ and all variables are quantified over $\{0, 1\}$. Without loss of generality we assume that the range of the pairing function $\langle \cdot, \cdot \rangle$ is disjoint from $0^*$.

**Proof** (of Claim 3.1)
Fix $A$. If (2) does not hold, then its negation holds, i.e,

$$(\exists \ell_y)_{y \in I_0}, (\forall r_y)_{y \in I_0}, (\exists \ell_y)_{y \in I_1}, (\forall r_y)_{y \in I_1},$$
$$\ldots, (\exists \ell_y)_{y \in I_{\log \beta(m)}}, (\forall r_y)_{y \in I_{\log \beta(m)}} : M_i^{A'}(0^m) \text{ rejects}. \qquad (4)$$

Switching the quantifiers $(\exists \ell_y)_{y \in I_j}$ and $(\forall r_y)_{y \in I_j}$ pairwise for every $0 \leqslant j \leqslant \log \beta(m)$ in (4) yields the weaker statement (3).

(Claim 3.1) □

Figure 1 describes the $i$-th stage in the construction of the set $A$. Note that the lexicographically first values in this algorithm always exist, so the construction works fine. We now argue that the resulting set $A$ satisfies the properties of Lemma 3.4.

1. The construction guarantees that by the end of stage $i$, $A(0^m) = 1 - M_i^{A \setminus \{0^m\}}(0^m)$ holds. Since $M_i$ on input $0^m$ cannot query $0^m$ (because $M_i$ is an autoreduction) nor any of the strings added during subsequent stages (because $M_i$ does not even have the time to write down any of these strings), $A(0^m) = 1 - M_i^A(0^m)$ holds for the final set $A$. So, $M_i$ is not an autoreduction of $A$. Since this is true of any autoreduction $M_i$, the set $A$ is not autoreducible.

2. During stage $i$, we encode $K|_I$ in the left region iff we do not put $0^m$ into $A$; otherwise we encode $K|_I$ in the right region. So, for any $y \in I$, $K(y) = A(\langle A(0^m), y \rangle)$. Therefore, $K \leqslant_{2-T}^{P} A$.

3. First note that $A$ only contains strings of the form $0^m$, and $\langle b, y \rangle$ for $b \in \Sigma$ and $y \in \Sigma^*$.

   The computation intensive steps during the $i$-th stage of the construction of $A$ are deciding $\text{QBF}_{\log \beta(m)}$-formulae of size $O(2^{m^{\beta(m)}})$ like (2) and deciding $K$ on inputs from $I$. All stages up to but not including $i$ can be done within time $O(2^m \cdot (t(m) + t'(m)))$. To decide $0^m$, we additionally have to check the validity of (2). To decide $\langle b, y \rangle$ with $b \in \Sigma$ and $y \in I_k$, we also have to do this, and we either have to evaluate $K(y)$ or run the part of stage $i$ corresponding to values of $j$ in Figure 1 up to and including $k$. The latter we can perform within time $O(2^{n^2} \cdot (t(n^2) + t'(n)))$, where $n = |y|$. A crucial point here is that $n^2$ upper bounds the length of the elements of $I_k$.

   All together, this yields the time bound claimed for $A$. The analysis of the space complexity is analogous.

> **if** formula (2) holds
>   **then for** $j = 0, \ldots, \log \beta(m)$
>       $(\ell_y)_{y \in I_j} \leftarrow (K(y))_{y \in I_j}$
>       $(r_y)_{y \in I_j} \leftarrow$ the lexicographically first value satisfying
>           $(\forall \ell_y)_{y \in I_{j+1}}, (\exists r_y)_{y \in I_{j+1}}, (\forall \ell_y)_{y \in I_{j+2}}, (\exists r_y)_{y \in I_{j+2}},$
>           $\ldots, (\forall \ell_y)_{y \in I_{\log \beta(m)}}, (\exists r_y)_{y \in I_{\log \beta(m)}} : M_i^{A'}(0^m)$ accepts,
>       where $A' = A \cup \{\langle 0, y \rangle \mid y \in I \text{ and } \ell_y = 1\} \cup \{\langle 1, y \rangle \mid y \in I \text{ and } r_y = 1\}$
>     **endfor**
>     $A \leftarrow A \cup \{\langle 0, y \rangle \mid y \in I \text{ and } \ell_y = 1\} \cup \{\langle 1, y \rangle \mid y \in I \text{ and } r_y = 1\}$
>   **else** { formula (3) holds }
>     **for** $j = 0, \ldots, \log \beta(m)$
>       $(r_y)_{y \in I_j} \leftarrow (K(y))_{y \in I_j}$
>       $(\ell_y)_{y \in I_j} \leftarrow$ the lexicographically first value satisfying
>           $(\forall r_y)_{y \in I_{j+1}}, (\exists \ell_y)_{y \in I_{j+1}}, (\forall r_y)_{y \in I_{j+2}}, (\exists \ell_y)_{y \in I_{j+2}},$
>           $\ldots, (\forall r_y)_{y \in I_{\log \beta(m)}}, (\exists \ell_y)_{y \in I_{\log \beta(m)}} : M_i^{A'}(0^m)$ accepts,
>       where $A' = A \cup \{\langle 0, y \rangle \mid y \in I \text{ and } \ell_y = 1\} \cup \{\langle 1, y \rangle \mid y \in I \text{ and } r_y = 1\}$
>     **endfor**
>     $A \leftarrow A \cup \{0^m\} \cup \{\langle 0, y \rangle \mid y \in I \text{ and } \ell_y = 1\} \cup \{\langle 1, y \rangle \mid y \in I \text{ and } r_y = 1\}$
> **endif**

Figure 1: Stage $i$ of the construction of the set $A$ in Lemma 3.4

Using the upper bound $2^{n^{\beta(n)}}$ for $\sigma'(n)$, the smallest standard complexity class to which Lemma 3.4 applies, seems to be EEXPSPACE. This results in Theorem 3.1.

**Proof** (of Theorem 3.1)
In Lemma 3.4, set $K$ a $\leqslant_m^P$-complete set for EEXPSPACE, and $\beta(n) = n$. □

In section 4.2, we will see that $\leqslant_{2-T}^P$ in the statement of Theorem 3.1 is optimal: Theorem 4.5 shows that Theorem 3.1 fails for $\leqslant_{2-tt}^P$.

We note that the proof of Theorem 3.1 carries through for $\leqslant_T^{\text{EXPSPACE}}$-reductions with polynomially bounded query lengths. This implies the strengthening given by Theorem 3.2.

## 3.2 Nonadaptive Autoreductions

In case of a nonadaptive autoreduction $M$ running in time $\tau(n)$, when diagonalizing against it on an input of length $n$, there are only $\tau(n)$ coding strings that can interfere, as opposed to $2^{\tau(n)}$ for autoreductions. This allows to reduce the complexity of the set constructed in Lemma 3.4 as follows:

**Lemma 3.5** *Fix a set $K$, and suppose we can decide it simultaneously in time $t(n)$ and space $\sigma(n)$. Let $\beta : \mathbb{N} \to (0, \infty)$ be a constructible monotone unbounded function, and suppose there is a deterministic Turing machine accepting TQBF that takes time $t'(n)$ and space $\sigma'(n)$ on QBF-formulae of size $n^{\beta(n)}$ with at most $\log \beta(n)$ alternations. Then there is a set $A$ such that:*

1. *$A$ is not nonadaptively autoreducible.*

2. *$K \leqslant_{3-tt}^P A$.*

3. *We can decide $A$ simultaneously in time $O(2^n \cdot (t(n^2) + t'(n)))$ and space $O(2^n \cdot (\sigma(n^2) + \sigma'(n)))$.*

**Proof** (of Lemma 3.5)
The construction of the set $A$ is the same as in Lemma 3.4 (see Figure 1) apart from the following differences:

- $M_1, M_2, \ldots$ now is a standard enumeration of nonadaptive autoreductions clocked such that $M_i$ runs in time $n^{\beta(i)}$ on inputs of length $n$.

- During stage $i \geqslant 1$ of the construction, $I$ denotes the set of all strings in $Q_{M_i}(0^m)$ with lengths in $[m, m^{\beta(m)} + 1) = [n_i, n_{i+1})$, and $I_j$ for $0 \leqslant j \leqslant \beta(m)$ denotes the set of strings in $Q_{M_i}(0^m)$ with lengths in $[m^{2^j}, \min(m^{2^{j+1}}, m^{\beta(m)} + 1))$.

- At the end of stage $i$, we additionally union $A$ with $\{\langle b, y \rangle \mid b \in \Sigma, y \in \Sigma^* \text{ with } m \leqslant |y| < m^{\beta(m)} + 1, y \notin I \text{ and } K(y) = 1\}$.

Adjusting time and space bounds appropriately, the proof that $A$ satisfies the 3 properties claimed, carries over. For the third one, note that (2) has become a $\text{QBF}_{\log \beta(n)}$-formula of length $O(n^{\beta(n)})$.
(Lemma 3.5) □

As a consequence, we can lower the space complexity in the equivalent of Theorem 3.1 from doubly exponential to singly exponential, yielding Theorem 3.3. In section 4.2, we will show we cannot reduce the number of queries from 3 to 2 in Theorem 3.3.

If we restrict the number of queries the nonadaptive autoreduction is allowed to make to some fixed polynomial, the proof technique of Theorem 3.3 also applies to EXP. In particular, we obtain:

**Theorem 3.6** *There is a $\leqslant^{\mathrm{P}}_{\mathrm{3-tt}}$-complete set for* EXP *that is not $\leqslant^{\mathrm{P}}_{\mathrm{btt}}$-autoreducible.*

## 4 Autoreducibility Results

For small complexity classes, all complete sets turn out to be autoreducible. Beigel and Feigenbaum [5] established this property of all levels of the polynomial-time hierarchy as well as of PSPACE, the largest class for which it was known to hold before our work. In this section, we will prove it for the $\Delta$-levels of the exponential-time hierarchy.

As to nonadaptive reductions, the question was even open for all levels of the polynomial-time hierarchy. We will show here that the $\leqslant^{\mathrm{P}}_{\mathrm{tt}}$-complete sets for the $\Delta$-levels of the polynomial-time hierarchy are nonadaptively autoreducible. For any complexity class containing EXP, we will prove that the $\leqslant^{\mathrm{P}}_{\mathrm{2-tt}}$-complete sets are $\leqslant^{\mathrm{P}}_{\mathrm{2-tt}}$-autoreducible.

Finally, we will also consider nonuniform and randomized autoreductions.

Throughout this section, we will assume without loss of generality an encoding $\gamma$ of a computation of a given oracle Turing machine $M$ on a given input $x$ with the following properties. $\gamma$ will be a marked concatenation of successive instantaneous descriptions of $M$, starting with the initial instantaneous description of $M$ on input $x$, such that:

- Given a pointer to a bit in $\gamma$, we can find out whether that bit represents the answer to an oracle query by probing a constant number of bits of $\gamma$.

- If it is the answer to an oracle query, the corresponding query is a substring of the prefix of $\gamma$ up to that point, and we can easily compute a pointer to the beginning of that substring without probing $\gamma$ any further.

- If it is not the answer to an oracle query, we can perform a local consistency check for that bit which only depends on a constant number of previous bit positions of $\gamma$ and the input $x$.

We call such an encoding a *valid computation* of $M$ on input $x$ iff the local consistency test for all the bit positions that do not correspond to oracle answers, are passed, and the other bits equal the oracle's answer to the corresponding query. Any other string we will call a *computation*.

### 4.1 Adaptive Autoreductions

We will first show that every $\leqslant^{\mathrm{P}}_{\mathrm{T}}$-complete set for EXP is autoreducible, and then generalize to all $\Delta$-levels of the exponential-time hierarchy.

**Theorem 4.1** *Every $\leqslant^{\mathrm{P}}_{\mathrm{T}}$-complete set for* EXP *is autoreducible.*

Here is the proof idea: For any of the standard deterministic complexity classes $\mathcal{C}$, we can decide each bit of the computation on a given input $x$ within $\mathcal{C}$. So, if $A$ is a $\leqslant^{\mathrm{P}}_{\mathrm{T}}$-complete set for $\mathcal{C}$, we can $\leqslant^{\mathrm{P}}_{\mathrm{T}}$-reduce this decision problem for $A$ to $A$. Now, consider the two (possibly invalid) computations

> **if** $R_\mu^{A\setminus\{x\}}(\langle x, 2^{p(|x|)}\rangle) = R_\mu^{A\cup\{x\}}(\langle x, 2^{p(|x|)}\rangle)$
>     **then accept** iff $R_\mu^{A\cup\{x\}}(\langle x, 2^{p(|x|)}\rangle) = 1$
>     **else** $i \leftarrow R_\sigma^{A\cup\{x\}}(x)$
>         **accept** iff local consistency test for $R_\mu^{A\setminus\{x\}}(\langle x, \cdot\rangle)$
>                         being the computation of $M$ on input $x$
>                         fails on the $i$-th bit
> **endif**

Figure 2: Autoreduction for the set $A$ of Theorem 4.1 on input $x$

we obtain by applying for every bit position the above reduction, answering all queries except for $x$ according to $A$, and assuming $x \in A$ for one computation, and $x \notin A$ for the other.

Note that the computation corresponding to the right assumption about $A(x)$, is certainly correct. So, if both computations yield the same answer (which we can efficiently check using $A$ without querying $x$), that answer is correct. If not, the other computation contains a mistake. We cannot check both computations entirely to see which one is right, but given a pointer to the first incorrect bit of the wrong computation, we can efficiently verify that it is mistaken by checking only a constant number of bits of that computation. The pointer is again computable within $\mathcal{C}$.

In case $\mathcal{C} \subseteq \text{EXP}$, using a $\leqslant_T^P$-reduction to $A$ and assuming $x \in A$ or $x \notin A$ as above, we can determine the pointer with oracle $A$ (but without querying $x$) in polynomial time, since the pointer's length is polynomially bounded.

We now fill out the details.

**Proof** (of Theorem 4.1)
Fix a $\leqslant_T^P$-complete set $A$ for EXP. Say $A$ is accepted by a Turing machine $M$ such that the computation of $M$ on an input of size $n$ has length $2^{p(n)}$ for some fixed polynomial $p$. Without loss of generality the last bit of the computation gives the final answer.

Let $\mu(\langle x, i\rangle)$ denote the $i$-th bit of the computation of $M$ on input $x$. We can compute $\mu$ in EXP, so there is an oracle Turing machine $R_\mu \leqslant_T^P$-reducing $\mu$ to $A$.

Let $\sigma(x)$ be the first $i$, $1 \leqslant i \leqslant 2^{p(|x|)}$, such that $R_\mu^{A\setminus\{x\}}(\langle x, i\rangle) \neq R_\mu^{A\cup\{x\}}(\langle x, i\rangle)$, provided it exists. Again, we can compute $\sigma$ in EXP, so there is a $\leqslant_T^P$-reduction $R_\sigma$ from $\sigma$ to $A$.

Consider the algorithm in Figure 2 for deciding $A$ on input $x$. The algorithm is a polynomial-time oracle Turing machine with oracle $A$ that does not query its own input $x$. We now argue that it correctly decides $A$ on input $x$. We distinguish between two cases:

**case** $R_\mu^{A\setminus\{x\}}(\langle x, 2^{p(|x|)}\rangle) = R_\mu^{A\cup\{x\}}(\langle x, 2^{p(|x|)}\rangle)$

Since at least one of the computations $R_\mu^{A\setminus\{x\}}(\langle x, \cdot\rangle)$ or $R_\mu^{A\cup\{x\}}(\langle x, \cdot\rangle)$ coincides with the actual computation of $M$ on input $x$, and the last bit of the computation equals the final decision, correctness follows.

**case** $R_\mu^{A\setminus\{x\}}(\langle x, 2^{p(|x|)}\rangle) \neq R_\mu^{A\cup\{x\}}(\langle x, 2^{p(|x|)}\rangle)$

If $x \in A$, $R_\mu^{A\setminus\{x\}}(\langle x, 2^{p(|x|)}\rangle) = 0$, so $R_\mu^{A\setminus\{x\}}(\langle x, \cdot\rangle)$ contains a mistake. Variable $i$ gets the correct value of the index of the first incorrect bit in this computation, so the local consistency test fails, and we accept $x$.

If $x \notin A$, $R_\mu^{A\setminus\{x\}}(\langle x, \cdot \rangle)$ is a valid computation, so no local consistency test fails, and we reject $x$.

(Theorem 4.1) □

The local checkability property of computations used in the proof of Theorem 4.1 does not relativize, because the oracle computation steps depend on the entire query, i.e., on a number of bits that is only limited by the resource bounds of the base machine, in this case exponentially many. We next show that Theorem 4.1 itself also does not relativize.

**Theorem 4.2** *Relative to some oracle,* EXP *has a* $\leqslant_{2-T}^P$*-complete set that is not autoreducible.*

**Proof**

Note that EXP has the following property:

**Property 4.1** *There is an oracle Turing machine $N$ running in* EXP *such that for any oracle $B$, the set accepted by $N^B$ is $\leqslant_m^P$-complete for* $\mathrm{EXP}^B$.

Without loss of generality, we assume that $N$ runs in time $2^n$. Let $K^B$ denote the set accepted by $N^B$.

We will construct an oracle $B$ and a set $A$ such that $A$ is $\leqslant_{2-T}^P$-complete for $\mathrm{EXP}^B$ and is not $\leqslant_T^{P^B}$-autoreducible.

The construction of $A$ is the same as in Lemma 3.4 (see Figure 1) with $\beta(n) = \log n$ and $K = K^B$, except for that the reductions $M_i$ now also have access to the oracle $B$.

We will encode in $B$ information about the construction of $A$ that reduces the complexity of $A$ relative to $B$, but do it high enough so as not to destroy the $\leqslant_{2-T}^P$-completeness of $A$ for $\mathrm{EXP}^B$ nor the diagonalizations against $\leqslant_T^{P^B}$-autoreductions.

We construct $B$ in stages along with $A$. We start with $B$ empty. Using the notation of Lemma 3.4, at the beginning of stage $i$, we add $0^{2^m}$ to $B$ iff property (2) does not hold, and at the end of sub-stage $j$, we union $B$ with

$$\{\langle 0^{2^{m^{2^{j+1}}}}, y\rangle \mid y \in I_j \text{ and } r(y) = 1\} \quad \text{if (2) holds at stage } i$$
$$\{\langle 0^{2^{m^{2^{j+1}}}}, y\rangle \mid y \in I_j \text{ and } \ell(y) = 1\} \quad \text{otherwise.}$$

Note that this does not affect the value of $K^B(y)$ for $|y| < m^{2^{j+1}}$, nor the computations of $M_i$ on inputs of size at most $m$ (for sufficiently large $i$ such that $m^{\log m} < 2^m$). It follows from the analysis in the proof of Lemma 3.4 that the set $A$ is $\leqslant_{2-T}^P$-hard for $\mathrm{EXP}^B$ and not $\leqslant_T^{P^B}$-autoreducible.

Regarding the complexity of deciding $A$ relative to $B$, note that the encoding in the oracle $B$ allows us to eliminate the need for evaluating $\mathrm{QBF}_{\log \beta(n)}$-formulae of size $2^{n^{\beta(n)}}$. Instead, we just query $B$ on easily constructed inputs of size $O(2^{n^2})$. Therefore, we can drop the terms corresponding to the $\mathrm{QBF}_{\log \beta(n)}$-formulae of size $2^{n^{\beta(n)}}$ in the complexity of $A$. Consequently, $A \in \mathrm{EXP}^B$.
(Theorem 4.2) □

Theorem 4.2 applies to any complexity class containing EXP that has Property 4.1, e.g., EXPSPACE, EEXP, EEXPSPACE, etc.

Sometimes, the structure of the oracle allows to get around the lack of local checkability of oracle queries. This is the case for oracles from the polynomial-time hierarchy, and leads to the following extension of Theorem 4.1:

**Theorem 4.3** *For any integer $k \geqslant 0$, every $\leqslant_{\mathrm{T}}^{\mathrm{P}}$-complete set for $\Delta_{k+1}^{\mathrm{EXP}}$ is autoreducible.*

The proof idea is as follows: Let $A$ be a $\leqslant_{\mathrm{T}}^{\mathrm{P}}$-complete set accepted by the deterministic oracle Turing machine $M$ with oracle $\mathrm{TQBF}_k$. First note that there is a polynomial-time Turing machine $N$ such that a query $q$ belongs to the oracle $\mathrm{TQBF}_k$ iff

$$\exists\, y_1, \forall\, y_2, \ldots, Q_k\, y_k : N(q, y_1, y_2, \ldots, y_k) \text{ accepts},\tag{5}$$

where the $y_i$'s are of size polynomial in $|q|$.

We consider the two purported computations of $M$ on input $x$ constructed in the proof of Theorem 4.1. One of them belongs to a party assuming $x \in A$, the other one to a party assuming $x \notin A$. The computation corresponding to the right assumption is correct; the other one might not be.

Now, suppose the computations differ, and we are given a pointer to the first bit position where they disagree, which turns out to be the answer to an oracle query $q$. Then we can have the two parties play the $k$-round game underlying (5): The party claiming $q \in \mathrm{TQBF}_k$ plays the existentially quantified $y_i$'s, the other one the universally quantified $y_i$'s. The players' strategies will consist of computing the game history so far, determining their optimal next move, $\leqslant_{\mathrm{T}}^{\mathrm{P}}$-reducing this computation to $A$, and finally producing the result of this reduction under their respective assumption about $A(x)$. This will guarantee that the party with the correct assumption plays optimally. Since this is also the one claiming the correct answer to the oracle query $q$, he will win the game, i.e., $N(q, y_1, y_2, \ldots, y_k)$ will equal his answer bit.

The only thing the autoreduction for $A$ has to do, is determine the value of $N(q, y_1, y_2, \ldots, y_k)$ in polynomial time using $A$ as an oracle but without querying $x$. It can do that along the lines of the base case algorithm given in Figure 2. If during this process, the local consistency test for $N$'s computation requires the knowledge of bits from the $y_i$'s, we compute these via the reduction defining the strategy of the corresponding player. The bits from $q$ we need, we can retrieve from the $M$-computations, since both computations are correct up to the point where they finished generating $q$. Once we know $N(q, y_1, y_2, \ldots, y_k)$, we can easily decide the correct assumption about $A(x)$.

The construction hinges on the hypothesis that we can $\leqslant_{\mathrm{T}}^{\mathrm{P}}$-reduce determining the player's moves to $A$. Computing these moves can become quite complex, though, because we have to recursively reconstruct the game history so far. The number of rounds $k$ being constant, seems crucial for keeping the complexity under control. The conference version of this paper [6] erroneously claimed the proof works for EXPSPACE, which can be thought of as alternating exponential time with an exponential number of alternations. Establishing Theorem 4.3 for EXPSPACE would actually separate NL from NP, as we will see in Section 5.

**Proof** (of Theorem 4.3)
Let $A$ be a $\leqslant_{\mathrm{T}}^{\mathrm{P}}$-complete set for $\Delta_{k+1}^{\mathrm{EXP}} = \mathrm{EXP}^{\Sigma_k^{\mathrm{P}}}$ accepted by the exponential-time oracle Turing machine $M$ with oracle $\mathrm{TQBF}_k$. Without loss of generality there is a polynomial $p$ and a polynomial-time Turing machine $N$ such that on inputs of size $n$, $M$ makes exactly $2^{p(n)}$ oracle queries, all of the form

$$\exists\, y_1 \in \Sigma^{2^{p(n)}}, \forall\, y_2 \in \Sigma^{2^{p(n)}}, \ldots, Q_k\, y_k \in \Sigma^{2^{p(n)}} : N(q, y_1, y_2, \ldots, y_k) \text{ accepts},\tag{6}$$

where $q$ has length $2^{p^2(n)}$. Moreover, the computations of $N$ in (6) each have length $2^{p^3(n)}$, and their last bit represents the answer; the same holds for the computations of $M$ on inputs of length $n$.

We first define a bunch of functions computable in $\Delta_{k+1}^{\text{EXP}}$. For each of them, say $\xi$, we fix an oracle Turing machine $R_\xi$ that $\leqslant_T^P$-reduces $\xi$ to $A$, and which the final autoreduction for $A$ will use. The proofs that we can compute these functions in $\Delta_{k+1}^{\text{EXP}}$ are straightforward.

Let $\mu(\langle x, i \rangle)$ denote the $i$-th bit of the computation of $M^{\text{TQBF}_k}$ on input $x$, and $\sigma(x)$ the first $i$ (if any) such that $R_\mu^{A \setminus \{x\}}(\langle x, i\rangle) \neq R_\mu^{A \cup \{x\}}(\langle x, i\rangle)$. The roles of $\mu$ and $\sigma$ are the same as in the proof of Theorem 4.1: We will use $R_\mu$ to figure out whether both possible answers for the oracle query "$x \in A$?" lead to the same final answer, and if not, use $R_\sigma$ to find a pointer $i$ to the first incorrect bit (in any) of the simulated computation getting the negative oracle answer $x \notin A$. If $i$ turns out not to point to an oracle query, we can proceed as in the proof of Theorem 4.1. Otherwise, we will make use of the following functions and associated reductions to $A$.

Let $1 \leqslant \ell \leqslant k$. The function $\eta_\ell$ takes an input $x$ such that the $i$-th bit of $R_\mu^{A \setminus \{x\}}(\langle x, \cdot \rangle)$ where $i = R_\sigma^{A \cup \{x\}}(x)$, is the answer to an oracle query, say (6). For $1 \leqslant m \leqslant \ell - 1$, let

$$y_m = \begin{cases} R_{\eta_m}^{A \cup \{x\}}(x) & \text{if } m \equiv R_\mu^{A \cup \{x\}}(\langle x, i \rangle) \bmod 2 \\ R_{\eta_m}^{A \setminus \{x\}}(x) & \text{otherwise.} \end{cases} \qquad (7)$$

We define $\eta_\ell(x)$ as the lexicographically least $y_\ell \in \Sigma^{2^{p(|x|)}}$ such that

$$\chi[Q_{\ell+1}\, y_{\ell+1}, Q_{\ell+2}\, y_{\ell+2}, \ldots, Q_k\, y_k : N(q, y_1, y_2, \ldots, y_k) \text{ accepts}] \equiv \ell \bmod 2;$$

if this value does not exist, we set $\eta_\ell(x) = 0^{2^{p(|x|)}}$.

In case $i$ points to the answer to an oracle query (6), the function $\eta_\ell$ and the reduction $R_{\eta_\ell}$ incorporate the moves during the $\ell$-th round of the game underlying (6). The function $\eta_\ell$ defines an optimal move during that round, provided it exists. The reduction $R_{\eta_\ell}$ together with the player's assumption about membership of $x$ to $A$, determines the actual move, which may differ from the one given by $\eta_\ell$ in case the player's assumption is incorrect. The condition on the right-hand side of (7) ensures that each opponent plays his rounds.

Finally, we define the functions $\nu$ and $\tau$, which have a similar job as the functions $\mu$ respectively $\sigma$, but for the computation of $N(q, y_1, y_2, \ldots, y_k)$ instead of the computation of $M_k^{\text{TQBF}}(x)$. More precisely, $\nu(\langle x, r \rangle)$ equals the $r$-th bit of the computation of $N(q, y_1, y_2, \ldots, y_k)$, where the $y_m$'s are defined by (7), and the bit with index $i = R_\sigma^{A \cup \{x\}}(x)$ in the computation $R_\mu^{A \setminus \{x\}}(\langle x, \cdot \rangle)$ is the answer to the oracle query (6). We define $\tau(x)$ to be the first $r$ (if any) for which $R_\nu^{A \setminus \{x\}}(\langle x, r \rangle) \neq R_\nu^{A \cup \{x\}}(\langle x, r \rangle)$, provided the bit with index $i = R_\sigma^{A \cup \{x\}}(x)$ in the computation $R_\mu^{A \setminus \{x\}}(\langle x, \cdot \rangle)$ is the answer to an oracle query.

Now we have these functions and corresponding reductions, we can describe an autoreduction for $A$. On input $x$, it works as described in Figure 3. We next argue that the algorithm correctly decides $A$ on input $x$. Checking the other properties required of an autoreduction for $A$ is straightforward.

We only consider the cases where $R_\mu^{A \setminus \{x\}}(\langle x, 2^{p^3(|x|)}\rangle) \neq R_\mu^{A \cup \{x\}}(\langle x, 2^{p^3(|x|)}\rangle)$ and $i$ points to the answer to an oracle query in $R_\mu^{A \setminus \{x\}}(\langle x, \cdot \rangle)$. We refer to the analysis in the proof of Theorem 4.1 for the remaining cases.

**case** $R_\nu^{A \setminus \{x\}}(\langle x, 2^{p^3(|x|)}\rangle) = R_\nu^{A \cup \{x\}}(\langle x, 2^{p^3(|x|)}\rangle)$

If $x \in A$, variable $i$ points to the first incorrect bit of $R_\mu^{A \setminus \{x\}}(\langle x, \cdot \rangle)$, which turns out to be the answer to an oracle query, say (6). Since the query is correctly described in $R_\mu^{A \setminus \{x\}}(\langle x, \cdot \rangle)$,

```
if R_μ^{A\{x}}(⟨x, 2^{p^3(|x|)}⟩) = R_μ^{A∪{x}}(⟨x, 2^{p^3(|x|)}⟩)
   then accept iff R_μ^{A∪{x}}(⟨x, 2^{p^3(|x|)}⟩) = 1
   else i ← R_σ^{A∪{x}}(x)
       if the i-th bit of R_μ^{A\{x}}(⟨x, ·⟩) is not the answer to an oracle query
           then accept iff local consistency test for R_μ^{A\{x}}(⟨x, ·⟩)
                               being the computation of M on input x
                               fails on the i-th bit
           else if R_ν^{A\{x}}(⟨x, 2^{p^3(|x|)}⟩) = R_ν^{A∪{x}}(⟨x, 2^{p^3(|x|)}⟩)
               then accept iff R_μ^{A\{x}}(⟨x, i⟩) ≠ R_ν^{A\{x}}(⟨x, 2^{p^3(|x|)}⟩)
               else r ← R_τ^{A∪{x}}(x)
                   accept iff local consistency test for (R_ν^{A\{x}}(⟨x, ·⟩)
                       being the computation of N(q, y_1, y_2, ..., y_k)
                       fails on the r-th bit,
                       where q denotes the query described in R_μ^{A\{x}}(⟨x, ·⟩)
                           to which the i-th bit in this computation is the answer
                       and
                       y_m = { R_{η_m}^{A∪{x}}(x)   if m ≡ R_μ^{A∪{x}}(⟨x, i⟩) mod 2
                             { R_{η_m}^{A\{x}}(x)   otherwise
               endif
           endif
       endif
```

Figure 3: Autoreduction for the set $A$ of Theorem 4.3 on input $x$

and the setting of the $y_m$'s forces the outcome of the game underlying (6) with input $q$ to the correct oracle answer $R_\mu^{A\cup\{x\}}(\langle x, i\rangle)$, $R_\nu^{A\cup\{x\}}(\langle x, 2^{p^3(|x|)}\rangle)$ gives the correct answer. Since $R_\mu^{A\setminus\{x\}}(\langle x, i\rangle)$ is incorrect,

$$R_\mu^{A\setminus\{x\}}(\langle x, i\rangle) \neq R_\nu^{A\cup\{x\}}(\langle x, 2^{p^3(|x|)}\rangle) = R_\nu^{A\setminus\{x\}}(\langle x, 2^{p^3(|x|)}\rangle),$$

and we accept $x$.

If $x \notin A$, both $R_\mu^{A\setminus\{x\}}(\langle x, i\rangle)$ and $R_\nu^{A\setminus\{x\}}(\langle x, 2^{p^3(|x|)}\rangle)$ give the correct answer to the oracle query $i$ points to in the computation $R_\mu^{A\setminus\{x\}}(\langle x, \cdot\rangle)$. So, we reject $x$.

**case** $R_\nu^{A\setminus\{x\}}(\langle x, 2^{p^3(|x|)}\rangle) \neq R_\nu^{A\cup\{x\}}(\langle x, 2^{p^3(|x|)}\rangle)$

If $x \in A$, $R_\nu^{A\setminus\{x\}}(\langle x, 2^{p^3(|x|)}\rangle)$ is incorrect, so $R_\nu^{A\setminus\{x\}}(\langle x, \cdot\rangle)$ has an error as a computation of $N(q, y_1, y_2, \ldots, y_k)$. Variable $r$ gets assigned the index of the first incorrect bit in this computation, so the local consistency check fails, and we accept $x$.

If $x \notin A$, $R_\nu^{A\setminus\{x\}}(\langle x,\cdot\rangle)$ is a valid computation of $N(q,y_1,y_2,\ldots,y_k)$, so every local consistency test is passed, and we reject $x$.

(Theorem 4.3) □

## 4.2 Nonadaptive Autoreductions

So far, we constructed autoreductions for $\leqslant_\mathrm{T}^\mathrm{P}$-complete sets $A$. On input $x$, we looked at the two candidate computations obtained by reducing to $A$, answering all oracle queries except for $x$ according to $A$, and answering query $x$ positively for one candidate, and negatively for the other. If the candidates disagreed, we tried to find out the right one, which there always was. We managed to get the idea to work for quite powerful sets $A$, e.g., EXP-complete sets, by exploiting the local checkability of computations. That allowed us to figure out the wrong computation without going through the entire computation ourselves: With help from $A$, we first computed a pointer to the first mistake in the wrong computation, and then verified it locally.

We cannot use this adaptive approach for constructing nonadaptive autoreductions. It seems like figuring out the wrong computation in a nonadaptive way, requires the autoreduction to perform the computation of the base machine itself, so the base machine has to run in polynomial time. Then checking the computation essentially boils down to verifying the oracle answers. Using the game characterization of the polynomial-time hierarchy along the same lines as in Theorem 4.3, we can do this for oracles from the polynomial-time hierarchy.

**Theorem 4.4** *For any integer $k \geqslant 0$, every $\leqslant_\mathrm{tt}^\mathrm{P}$-complete set for $\Delta_{k+1}^\mathrm{P}$ is nonadaptively autoreducible.*

Parallel to the adaptive case, an earlier version of this paper [6] stated Theorem 4.4 for unbounded $k$, i.e., for PSPACE. However, we only get the proof to work for constant $k$. In Section 5, we will see that proving Theorem 4.4 for PSPACE would separate NL from NP.

The only additional difficulty in the proof is that in the nonadaptive setting, we do not know which player has to perform the even rounds, and which one the odd rounds in the $k$-round game underlying a query like (5). But we can just have them play both scenarios, and afterwards figure out the relevant run.

**Proof** (of Theorem 4.4)
Let $A$ be a $\leqslant_\mathrm{tt}^\mathrm{P}$-complete set for $\Delta_{k+1}^\mathrm{P} = \mathrm{P}^{\Sigma_k^\mathrm{P}}$ accepted by the polynomial-time oracle Turing machine $M$ with oracle $\mathrm{TQBF}_k$. Without loss of generality there is a polynomial $p$ and a polynomial-time Turing machine $N$ such that on inputs of size $n$, $M$ makes exactly $p(n)$ oracle queries $q$, all of the form

$$\exists\, y_1 \in \Sigma^{p(n)},\, \forall\, y_2 \in \Sigma^{p(n)},\ldots, Q_k\, y_k \in \Sigma^{p(n)} : N(q,y_1,y_2,\ldots,y_k) \text{ accepts}, \qquad (8)$$

where $q$ has length $p^2(n)$. Let $q(x,i)$ denote the $i$-th oracle query of $M^{\mathrm{TQBF}_k}$ on input $x$. Note that $q \in \mathrm{FP}^{\Sigma_k^\mathrm{P}}$.

Let $Q = \{\langle x,i\rangle \mid q(x,i) \in \mathrm{TQBF}_k\}$. The set $Q$ belongs to $\Delta_{k+1}^\mathrm{P}$, so there is a $\leqslant_\mathrm{tt}^\mathrm{P}$-reduction $R_Q$ from $Q$ to $A$.

If for a given input $x$, $R_Q^{A\cup\{x\}}$ and $R_Q^{A\setminus\{x\}}$ agree on $\langle x,j\rangle$ for every $1 \leqslant j \leqslant p(|x|)$, we are home: We can simulate the base machine $M$ using $R_Q^{A\cup\{x\}}(\langle x,j\rangle)$ as the answer to the $j$-th oracle query.

Otherwise, we will make use of the following functions $\eta_1, \eta_2, \ldots, \eta_k$ computable in $\Delta_{k+1}^{\mathrm{P}}$, and corresponding oracle Turing machines $R_{\eta_1}, R_{\eta_2}, \ldots, R_{\eta_k}$ defining $\leqslant_{\mathrm{tt}}^{\mathrm{P}}$-reductions to $A$: Let $1 \leqslant \ell \leqslant k$. The function $\eta_\ell$ is defined for inputs $x$ such that there is a smallest $1 \leqslant i \leqslant p(|x|)$ for which $R_Q^{A\setminus\{x\}}(\langle x,i\rangle) \neq R_Q^{A\cup\{x\}}(\langle x,i\rangle)$. For $1 \leqslant m \leqslant \ell - 1$, let

$$y_m = \begin{cases} R_{\eta_m}^{A\cup\{x\}}(x) & \text{if } m \equiv R_Q^{A\cup\{x\}}(\langle x,i\rangle) \bmod 2 \\ R_{\eta_m}^{A\setminus\{x\}}(x) & \text{otherwise.} \end{cases} \quad (9)$$

We define $\eta_\ell(x)$ as the lexicographically least $y_\ell \in \Sigma^{p(|x|)}$ such that

$$\chi[Q_{\ell+1} y_{\ell+1}, Q_{\ell+2} y_{\ell+2}, \ldots, Q_k y_k : N(q(x,i), y_1, y_2, \ldots, y_k) \text{ accepts}] \equiv \ell \bmod 2;$$

we set $\eta_\ell(x) = 0^{p(|x|)}$ if such string does not exist.

The intuitive meaning of the functions $\eta_\ell$ and the reductions $R_{\eta_\ell}$ is similar to in the proof of Theorem 4.3: They capture the moves during the $\ell$-th round of the game underlying (8) for $q = q(x,i)$. The function $\eta_\ell$ encapsulates an optimal move during round $\ell$ if it exists, and the reduction $R_{\eta_\ell}$ under the player's assumption regarding membership of $x$ to $A$, produces the actual move in that round. The condition on the right-hand side of (9) guarantees the correct alternation of rounds.

Consider the algorithm in Figure 4. Note that the only queries to $A$ the algorithm in Figure 4

---

**if** $R_Q^{A\setminus\{x\}}(\langle x,j\rangle) = R_Q^{A\cup\{x\}}(\langle x,j\rangle)$ for every $1 \leqslant j \leqslant p(|x|)$
  **then accept** iff $M$ accepts $x$ when the $j$-th oracle query is answered $R_Q^{A\cup\{x\}}(\langle x,j\rangle)$
  **else** $i \leftarrow$ first $j$ such that $R_Q^{A\setminus\{x\}}(\langle x,j\rangle) \neq R_Q^{A\cup\{x\}}(\langle x,j\rangle)$
    **accept** iff $N(q, y_1, y_2, \ldots, y_k) = R_Q^{A\cup\{x\}}(\langle x,i\rangle)$
      where $q$ denotes the $i$-th query of $M$ on input $x$
        when the answer to the $j$-th oracle query is given by $R_Q^{A\cup\{x\}}(\langle x,j\rangle)$
      and
$$y_m = \begin{cases} R_{\eta_m}^{A\cup\{x\}}(x) & \text{if } m \equiv R_Q^{A\cup\{x\}}(\langle x,i\rangle) \bmod 2 \\ R_{\eta_m}^{A\setminus\{x\}}(x) & \text{otherwise} \end{cases}$$
**endif**

Figure 4: Nonadaptive autoreduction for the set $A$ of Theorem 4.4 on input $x$

needs to make, are the queries of $R_Q$ different from $x$ on inputs $\langle x,j\rangle$ for $1 \leqslant j \leqslant p(|x|)$, and the queries of $R_{\eta_m}$ different from $x$ on input $x$ for $1 \leqslant m \leqslant k$. Since $R_Q$ and the $R_{\eta_m}$'s are nonadaptive, it follows that Figure 4 describes a $\leqslant_{\mathrm{tt}}^{\mathrm{P}}$-reduction to $A$ that does not query its own input. A similar but simplified argument as in the proof of Theorem 4.3 shows that it accepts $A$. So, $A$ is nonadaptively autoreducible. (Theorem 4.4) □

Next, we consider more restricted reductions. Using a different technique, we show:

**Theorem 4.5** *For any complexity class $\mathcal{C}$, every $\leqslant_{2-\mathrm{tt}}^{\mathrm{P}}$-complete set for $\mathcal{C}$ is $\leqslant_{2-\mathrm{tt}}^{\mathrm{P}}$-autoreducible, provided $\mathcal{C}$ is closed under exponential-time reductions that only ask one query which is smaller in length.*

In particular, Theorem 4.5 applies to $\mathcal{C} = $ EXP, EXPSPACE, and EEXPSPACE. In view of Theorems 3.1 and 3.3, this implies that Theorems 3.1, 3.3, and 4.5 are optimal.

The proof exploits the ability of EXP to simulate all polynomial-time reductions to construct an auxiliary set $D$ within $\mathcal{C}$ such that any $\leqslant^{\mathrm{P}}_{\text{2-tt}}$-reductions of $D$ to some fixed complete set $A$ has a property that induces an autoreduction on $A$.

**Proof** (of Theorem 4.5)
Let $M_1, M_2, \ldots$ be a standard enumeration of $\leqslant^{\mathrm{P}}_{\text{2-tt}}$-reductions such that $M_i$ runs in time $n^i$ on inputs of size $n$. Let $A$ be a $\leqslant^{\mathrm{P}}_{\text{2-tt}}$-complete set for $\mathcal{C}$.

Consider the set $D$ that only contains strings of the form $\langle 0^i, x \rangle$ for $i \in \mathbb{N}$ and $x \in \Sigma^*$, and is decided by the algorithm of Figure 5 on such an input. Except for deciding $A(x)$, the algorithm

---

**case** truth-table of $M_i$ on input $\langle 0^i, x \rangle$ with the truth-value of query $x$ set to $A(x)$

    constant:

        **accept** iff $M_i^A$ rejects $\langle 0^i, x \rangle$

    of the form "$y \notin A$":

        **accept** iff $x \notin A$

    otherwise:

        **accept** iff $x \in A$

**endcase**

---

Figure 5: Algorithm for the set $D$ of Theorem 4.5 on input $\langle 0^i, x \rangle$.

runs in exponential time. Therefore, under the given conditions on $\mathcal{C}$, $D \in \mathcal{C}$, so there is a $\leqslant^{\mathrm{P}}_{\text{2-tt}}$-reduction $M_j$ from $D$ to $A$.

The construction of $D$ diagonalizes against every $\leqslant^{\mathrm{P}}_{\text{2-tt}}$-reduction $M_i$ of $D$ to $A$ whose truth-table on input $\langle 0^i, x \rangle$ would become constant once we filled in the membership bit for $x$. Therefore, for every input $x$, one of the following cases holds:

- The reduced truth-table is of the form "$y \in A$" with $y \neq x$.
  Then $y \in A \Leftrightarrow M_j$ accepts $\langle 0^j, x \rangle \Leftrightarrow x \in A$.

- The reduced truth-table is of the form "$y \notin A$" with $y \neq x$.
  Then $y \notin A \Leftrightarrow M_j$ accepts $\langle 0^j, x \rangle \Leftrightarrow x \notin A$.

- The truth-table depends on the membership to $A$ of 2 strings different from $x$.
  Then $M_j^A$ does not query $x$ on input $\langle 0^j, x \rangle$, and accepts iff $x \in A$.

The above analysis shows that the algorithm of Figure 6 describes a $\leqslant^{\mathrm{P}}_{\text{2-tt}}$-reduction of $A$.
(Theorem 4.5) □

### 4.3 Probabilistic and Nonuniform Autoreductions

The previous results in this section trivially imply that the $\leqslant^{\mathrm{P}}_{\mathrm{T}}$-complete sets for the $\Delta$-levels of the exponential-time hierarchy are probabilistically autoreducible, and the $\leqslant^{\mathrm{P}}_{\mathrm{tt}}$-complete sets for the $\Delta$-

```
if |Q_{M_j}(⟨0^j, x⟩) \ {x}| = 2
    then accept iff M_j^A accepts ⟨0^j, x⟩
    else { |Q_{M_j}(⟨0^j, x⟩) \ {x}| = 1 }
        y ← unique element of Q_{M_j}(⟨0^j, x⟩) \ {x}
        accept iff y ∈ A
endif
```

Figure 6: Autoreduction constructed in the proof of Theorem 4.5

levels of polynomial-time hierarchy are probabilistically nonadaptively autoreducible. Randomness allows us the prove more in the nonadaptive case.

First, we can establish Theorem 4.4 for EXP:

**Theorem 4.6** *Let $f$ be a constructible function. Every $\leqslant^{\mathrm{P}}_{f(n)-\mathrm{tt}}$-complete set for EXP is probabilistically $\leqslant^{\mathrm{P}}_{O(f(n))-\mathrm{tt}}$-autoreducible. In particular, every $\leqslant^{\mathrm{P}}_{\mathrm{tt}}$-complete set for EXP is probabilistically nonadaptively autoreducible.*

**Proof** (of Theorem 4.6)
Let $A$ be a $\leqslant^{\mathrm{P}}_{f(n)-\mathrm{tt}}$-complete set for EXP. We will apply the PCP Theorem for EXP [2] to $A$.

**Lemma 4.7 ([2])** *There is a constant $k$ such that for any set $A \in$ EXP, there is a polynomial-time Turing machine $V$ and a polynomial $p$ such that for any input $x$:*

- *If $x \in A$, then there exists a proof oracle $\pi$ such that*

$$\Pr_{|r|=p(|x|)}[V^\pi(x,r) \ accepts\,] = 1. \qquad (10)$$

- *If $x \notin A$, then for any proof oracle $\pi$*

$$\Pr_{|r|=p(|x|)}[V^\pi(x,r) \ accepts\,] \leqslant \frac{1}{3}.$$

*Moreover, $V$ never makes more than $k$ proof oracle queries, and there is a proof oracle $\tilde{\pi} \in$ EXP independent of $x$ such that (10) holds for $\pi = \tilde{\pi}$ in case $x \in A$.*

Translating Lemma 4.7 into our terminology, we obtain:

**Lemma 4.8** *There is a constant $k$ such that for any set $A \in$ EXP, there is a probabilistic $\leqslant^{\mathrm{P}}_{k-\mathrm{tt}}$-reduction $N$, and a set $B \in$ EXP such that for any input $x$:*

- *If $x \in A$, then $N^B(x)$ always accepts.*

- *If $x \notin A$, then for any oracle $C$, $N^C(x)$ accepts with probability at most $\frac{1}{3}$.*

Let $R$ be a $\leqslant^{\mathrm{P}}_{f(n)-\mathrm{tt}}$-reduction of $B$ to $A$, and consider the probabilistic reduction $M^A$ that on input $x$, runs $N$ on input $x$ with oracle $R^{A \cup \{x\}}$. $M^A$ is a probabilistic $\leqslant^{\mathrm{P}}_{k \cdot f(n)-\mathrm{tt}}$-reduction to $A$ that never queries its own input. The following shows it defines a reduction from $A$:

- If $x \in A$, $R^{A \cup \{x\}} = R^A = B$, so $M^A(x) = N^B(x)$ always accepts.

- If $x \notin A$, then for $C = R^{A \cup \{x\}}$, $M^A(x) = N^C(x)$ accepts with probability at most $\frac{1}{3}$.

(Theorem 4.6) □

Note that Theorem 4.6 makes it plausible why we did not manage to scale down Theorem 3.2 by one exponent to EXPSPACE in the nonadaptive setting, as we were able to do for our other results in Section 3 when going from the adaptive to the nonadaptive case: This would separate EXP from EXPSPACE.

We suggest the extension of Theorem 4.6 to the $\Delta$-levels of the exponential-time hierarchy as an interesting problem for further research.

Second, Theorem 4.4 also holds for NP:

**Theorem 4.9** *All $\leqslant_{\mathrm{tt}}^{\mathrm{P}}$-complete sets for* NP *are probabilistically nonadaptively autoreducible.*

**Proof** (of Theorem 4.9)
Fix a $\leqslant_{\mathrm{tt}}^{\mathrm{P}}$-complete set $A$ for NP. Let $R_A$ denote a length nondecreasing $\leqslant_{\mathrm{m}}^{\mathrm{P}}$-reduction of $A$ to SAT.

Define the set

$$W = \{\langle \phi, 0^i \rangle \mid \phi \text{ is a Boolean formula with, say } m \text{ variables and } \exists a \in \Sigma^m : [\phi(a) \text{ and } a_i = 1]\}.$$

Since $W \in \mathrm{NP}$, there is a $\leqslant_{\mathrm{tt}}^{\mathrm{P}}$-reduction $R_W$ from $W$ to $A$.

We will use the following probabilistic algorithm by Valiant and Vazirani [17]:

**Lemma 4.10 ([17])** *There exists a polynomial-time probabilistic Turing machine $N$ that on input a Boolean formula $\varphi$ with $n$ variables, outputs another quantifier free Boolean formula $\phi = N(\varphi)$ such that:*

- *If $\varphi$ is satisfiable, then with probability at least $\frac{1}{4n}$, $\phi$ has a unique satisfying assignment.*

- *If $\varphi$ is not satisfiable, then $\phi$ is never satisfiable.*

Now consider the following algorithm for $A$: On input $x$, run $N$ on input $R_A(x)$, yielding a Boolean formula $\phi$ with, say $m$ variables, and it accepts iff

$$\phi(R_W^{A \cup \{x\}}(\langle \phi, 0 \rangle), R_W^{A \cup \{x\}}(\langle \phi, 00 \rangle), \ldots, R_W^{A \cup \{x\}}(\langle \phi, 0^i \rangle), \ldots, R_W^{A \cup \{x\}}(\langle \phi, 0^m \rangle))$$

evaluates to true. Note that this algorithm describes a probabilistic $\leqslant_{\mathrm{tt}}^{\mathrm{P}}$-reduction to $A$ that never queries its own input. Moreover:

- If $x \in A$, then with probability at least $\frac{1}{4|x|}$, the Valiant-Vazirani algorithm $N$ produces a Boolean formula $\phi$ with a unique satisfying assignment $\tilde{a}_\phi$. In that case, the assignment we use $(R_W^{A \cup \{x\}}(\langle \phi, 0 \rangle), R_W^{A \cup \{x\}}(\langle \phi, 00 \rangle), \ldots, R_W^{A \cup \{x\}}(\langle \phi, 0^i \rangle), \ldots, R_W^{A \cup \{x\}}(\langle \phi, 0^m \rangle))$ equals $\tilde{a}_\phi$, and we accept $x$.

- If $x \notin A$, any Boolean formula $\phi$ which $N$ produces has no satisfying assignment, so we always reject $x$.

Executing $\Theta(n)$ independent runs of this algorithm, and accepting iff any of them accepts, yields a probabilistic nonadaptive autoreduction for $A$. (Theorem 4.9) □

So, for probabilistic autoreductions, we get similar results as for deterministic ones: Low end complexity classes turn out to have the property that their complete sets are autoreducible, whereas high end complexity classes do not. As we will see in more detail in the next section, this structural difference yields separations.

If we allow nonuniformity, the situation changes dramatically. Since probabilistic autoreducibility implies nonuniform autoreducibility [5], all our positive results for small complexity classes carry over to the nonuniform setting. But, as we will see next, the negative results do not, because also the complete sets for large complexity classes become autoreducible, both in the adaptive and in the nonadaptive case. So, uniformity is crucial for separating complexity classes using autoreducibility, and the Razborov-Rudich result [14] does not apply.

Feigenbaum and Fortnow [7] define the following concept of #P-robustness, of which we also consider the nonadaptive variant.

**Definition 4.1** *A set $A$ is* #P-*robust if* $\#P^A \subseteq FP^A$; *$A$ is* nonadaptively #P-*robust if* $\#P^A_{tt} \subseteq FP^A_{tt}$.

Nonadaptive #P-robustness implies #P-robustness. For the usual deterministic and nondeterministic complexity classes containing PSPACE, all $\leqslant^P_T$-complete sets are #P-robust. For the deterministic classes containing PSPACE, it is also true that the $\leqslant^P_{tt}$-complete sets are nonadaptively #P-robust.

The following connection with nonuniform autoreducibility holds:

**Theorem 4.11** *All* #P-*robust sets are nonuniformly autoreducible. All nonadaptively* #P-*robust sets are nonuniformly nonadaptively autoreducible.*

**Proof**
Feigenbaum and Fortnow [7] show that every #P-robust language is random-self-reducible. Beigel and Feigenbaum [5] prove that every random-self-reducible set is nonuniformly autoreducible (or "weakly coherent" as they call it). Their proofs carry over to the nonadaptive setting. □

It follows that the $\leqslant^P_{tt}$-complete sets for the usual deterministic complexity classes containing PSPACE are all nonuniformly nonadaptively autoreducible. The same holds for adaptive reductions, in which case the property is also true of nondeterministic complexity classes containing PSPACE. In particular, we get the following:

**Corollary 4.12** *All $\leqslant^P_T$-complete sets for* NEXP, EXPSPACE, EEXP, NEEXP, EEXPSPACE, *... are nonuniformly autoreducible. All $\leqslant^P_{tt}$-complete sets for* PSPACE, EXP, EXPSPACE, *... are nonuniformly nonadaptively autoreducible.*

## 5 Separation Results

In this section, we will see how we can use the structural property of all complete sets being autoreducible to separate complexity classes. Based on the results of Sections 3 and 4, we only get separations that were already known: EXPH $\neq$ EEXPSPACE (by Theorems 4.3 and 3.1),

EXP ≠ EEXPSPACE (by Theorems 4.6 and 3.2), and PH ≠ EXPSPACE (by Theorems 4.4 and 3.3, and also by scaling down EXPH ≠ EEXPSPACE). However, settling the question for certain other classes, would yield impressive new separations.

We summarize the implications in Figure 7.

**Theorem 5.1** *In Figure 7, a positive answer to a question from the first column, implies the separation in the second column, and a negative answer, the separation in the third column.*

| question | yes | no |
| --- | --- | --- |
| Are all $\leqslant_{\mathrm{T}}^{\mathrm{P}}$-complete sets for EXPSPACE autoreducible? | NL ≠ NP | PH ≠ PSPACE |
| Are all $\leqslant_{\mathrm{T}}^{\mathrm{P}}$-complete sets for EEXP autoreducible? | NL ≠ NP <br> P ≠ PSPACE | PH ≠ EXP |
| Are all $\leqslant_{\mathrm{tt}}^{\mathrm{P}}$-complete sets for PSPACE $\leqslant_{\mathrm{tt}}^{\mathrm{P}}$-autoreducible? | NL ≠ NP | PH ≠ PSPACE |
| Are all $\leqslant_{\mathrm{tt}}^{\mathrm{P}}$-complete sets for EXP $\leqslant_{\mathrm{tt}}^{\mathrm{P}}$-autoreducible? | NL ≠ NP <br> P ≠ PSPACE | PH ≠ EXP |
| Are all $\leqslant_{\mathrm{tt}}^{\mathrm{P}}$-complete sets for EXPSPACE probabilistically $\leqslant_{\mathrm{tt}}^{\mathrm{P}}$-autoreducible? | NL ≠ NP | P ≠ PSPACE |

Figure 7: Separation results using autoreducibility

Most of the entries in Figure 7 follow directly from the results of the previous sections. In order to finish the table, we use the next lemma:

**Lemma 5.2** *If* NP = NL, *we can decide the validity of QBF-formulae of size $t$ and with $\alpha$ alternations on a deterministic Turing machine $M_1$ in time $t^{O(c^\alpha)}$ and on a nondeterministic Turing machine $M_2$ in space $O(c^\alpha \log t)$, for some constant $c$.*

**Proof** (of Lemma 5.2)
Since coNP = NP, by Cook's Theorem we can transform in polynomial time a $\Pi_1$-formula with free variables into an equivalent $\Sigma_1$-formula with the same free variables, and vice versa. Since NP = P, we can decide the validity of $\Sigma_1$-formulae in polynomial-time. Say both the transformation algorithm $T$ and the satisfiability algorithm $S$ run in time $n^c$ for some constant $c$.

Let $\phi$ be a QBF-formula of size $t$ with $\alpha$ alternations. Consider the following algorithm for deciding $\phi$: Repeatedly apply the transformation $T$ to the largest suffix that constitutes a $\Sigma_1$- or $\Pi_1$-formula until the whole formula becomes $\Sigma_1$, and then run $S$ on it.

This algorithm correctly decides the truth of $\phi$. Since the number of alternations decreases by one during every iteration, it makes at most $\alpha$ calls to $T$, each time at most raising the length of the formula to the power $c$. It follows that the algorithm runs in time $t^{O(c^\alpha)}$.

Moreover, since P = NL, a padding argument shows that DTIME[$\tau$] ⊆ NSPACE[log $\tau$] for any time constructible function $\tau$. Therefore the result holds. (Lemma 5.2) □

This allows us to improve Theorems 3.2 and 3.3 as follows under the hypothesis NP = NL:

**Theorem 5.3** *If* NP = NL, *there is a* $\leqslant^P_{2-T}$*-complete set for* EXPSPACE *that is not probabilistically autoreducible. The same holds for* EEXP *instead of* EXPSPACE.

**Proof**
Combine Lemma 5.2 with the probabilistic extension of Lemma 3.4 used in the proof of Theorem 3.2.

We have studied the question whether all complete sets are autoreducible for various complexity

**Theorem 5.4** *If* NP = NL, *there is a* $\leqslant^P_{3-tt}$*-complete set for* PSPACE *that is not nonadaptively autoreducible. The same holds for* EXP *instead of* PSPACE.

**Proof**
Combine Lemma 5.2 with Lemma 3.5. □

Now, we have all ingredients for establishing Figure 7:

**Proof** (of Theorem 5.1)
The NL $\neq$ NP implications in the "yes"-column of Figure 7 immediately follow from Theorems 5.3 and 5.4 by contraposition.

By Theorem 3.1, a positive answer to the 2nd question in Figure 7 would yield EEXP $\neq$ EEXPSPACE, and by Theorem 3.3, a positive answer to the 4th question would imply EXP $\neq$ EXPSPACE. By padding, both translate down to P $\neq$ PSPACE.

Similarly, by Theorem 4.3, a negative answer to the 2nd question would imply EXPH $\neq$ EEXP, which pads down to PH $\neq$ EXP. A negative answer to the 4th question would yield PH $\neq$ EXP directly by Theorem 4.4. By the same token, a negative answer to the 1st question results in EXPH $\neq$ EXPSPACE and PH $\neq$ PSPACE, and a negative answer to the 3rd question in PH $\neq$ PSPACE. By Theorem 4.6, a negative answer to the last question implies EXP $\neq$ EXPSPACE and P $\neq$ PSPACE. □

We note that we can tighten all of the separations in Figure 7 a bit, because we can apply Lemmata 3.4 and 3.5 to smaller classes than in Theorems 3.1 respectively 3.3. One improvement along these lines that might warrant attention is that we can replace "NL $\neq$ NP" in Figure 7 by "coNP $\not\subseteq$ NP $\cap$ NSPACE[$\log^{O(1)} n$]." This is because that condition suffices for Theorems 5.3 and 5.4, since we can strengthen Lemma 5.2 as follows:

**Lemma 5.5** *If* coNP $\subseteq$ NP $\cap$ NSPACE[$\log^{O(1)} n$], *we can decide the validity of QBF-formulae of size* $t$ *and with* $\alpha$ *alternations on a deterministic Turing machine* $M_1$ *in time* $t^{O(c^\alpha)}$ *and on a nondeterministic Turing machine* $M_2$ *in space* $O(d^\alpha \log^d t)$, *for some constants* $c$ *and* $d$.

# 6 Conclusion

classes and various reducibilities. We obtained a positive answer for lower complexity classes in Section 4, and a negative one for higher complexity classes in Section 3. This way, we separated these lower complexity classes from these higher ones by highlighting a structural difference. The resulting separations were not new, but we argued in Section 5 that settling the very same question for intermediate complexity classes, would provide major new separations.

We believe that refinements to our techniques may lead to them, and would like to end with a few words about some thoughts in that direction.

One does not have to look at complete sets only. Let $\mathcal{C}_1 \subseteq \mathcal{C}_2$. Suppose we know that all complete sets for $\mathcal{C}_2$ are autoreducible. Then it suffices to construct, e.g., along the lines of Lemma 3.4, a hard set for $\mathcal{C}_1$ that is not autoreducible, in order to separate $\mathcal{C}_1$ from $\mathcal{C}_2$.

As we mentioned at the end of Section 5, we can improve Theorem 3.1 a bit by applying Lemma 3.4 to smaller space-bounded classes than EEXPSPACE. We can not hope to gain much, though, since the coding in the proof of Lemma 3.4 seems to be $\text{DSPACE}[2^{n^{\beta(n)}}]$-complete because of the $\text{QBF}_{\log \beta(n)}$-formulae of size $2^{n^{\beta(n)}}$ involved for inputs of size $n$. The same holds for Theorem 3.3 and Lemma 3.5.

Generalizations of autoreducibility may allow us to push things further. For example, one could look at $k(n)$-autoreducibility where $k(n)$ bits of the set remain unknown to the querying machine. Theorem 4.3 goes through for $k(n) \in O(\log n)$. Perhaps one can exploit this leeway in the coding of Lemma 3.4 and narrow the gap between the positive and negative results. As discussed in Section 5, that would yield interesting separations.

Finally, one may want to look at other properties than autoreducibility to realize Post's Program in complexity theory. Perhaps another concept from computability theory or a more artificial property can be used to separate complexity classes.

## Acknowledgments

We would like to thank Manindra Agrawal and Ashish Naik for very helpful discussions. We are also grateful to Carsten Lund and Muli Safra for their help regarding the PCP Theorem.